\documentstyle[12pt]{article}
\begin{document}
\newtheorem{guess}{Proposition }[section]
\newtheorem{theorem}[guess]{Theorem}
\newtheorem{lemma}[guess]{Lemma}
\newtheorem{corollary}[guess]{Corollary}
\def \ja {\vrule height 3mm width 3mm}

\centerline{\bf\Large Holomorphic bundles on   ${\cal O}(-k)$ are algebraic}
\centerline { Elizabeth Gasparim}
\vspace{5 mm}

\begin{abstract}

We show that holomorphic bundles on  ${\cal O}(-k)$ for $k > 0$
 are algebraic. We also show
 holomorphic  bundles on O(-1) are trivial outside the zero section.

\end{abstract}

\section{ Preliminaires }

The line  bundle on ${\bf P}^1$ given 
by transition function $z^{k}$ is usually denoted ${\cal O} (-k).$
Since we will be studying bundles
over this space, we will denote $ {\cal O} (-k)$ by $M_k$
when we want to view this space as the base of a bundle. 
We give $M_k$ the following charts
 $M_k = U \cup V$, for
$$U =   {\bf C}^2 =\{(z,u)\}$$
$$V =   {\bf C}^2 = \{(\xi,v)\}$$
$$U \cap V = ({\bf  C} - \{0\})  \times   {\bf C}$$
with change of coordinates
$$ (\xi,v) = (z^{-1 },z^ku).$$

 Since  $H^1( {\cal O}(-k), {\cal O}) = 0,$
using the exponential sheaf sequence 
it follows that $Pic({\cal O}(-k)) =  {\bf Z},$
and 
holomorphic line bundles on $M_k$ 
are classified by their Chern classes.  Therefore
it is clear that holomorphic line bundles 
over $M_k$ are algebraic.
We will denote by  ${\cal O}^!(j) $ the line bundle on 
$M_k$ given by transition funcion $z^{-j}.$

If $E$ is a rank $n$ bundle over $M_k,$ 
then over the zero section (which is a 
 ${\bf P}^1$) $E$
splits as a sum of line bundles by Grothendieck's 
theorem.
Denoting the zero section  by $\ell$
it follows that for some integers $j_i$
uniquely determined up to order
$E_{\ell} \simeq  \oplus_{i = 1}^n {\cal O}(j_i).   $
We will show that such $E$ is an algebraic extension 
of the line bundles $ {\cal O}^!(j_i).$

\section{ Bundles on
${\cal O}(-k)$ are algebraic}

\begin{lemma}: Holomorphic bundles on $M_k$ are extensions 
of line bundles.
\end{lemma}

\noindent{\bf Proof}:
We give the proof for rank two for simplicity. The case for rank $n$
is proved by induction on $n$ using similar calculations.
Suppose rank $E$ = 2 and 
$E_{\ell} \simeq {\cal O}(-j_1) \oplus  {\cal O}(-j_2)$ 
 which we may assume to satisfy
$j_1 \ge j_2.$
A  transition matrix for $E$
from $U$ to $V$ therefore takes the form
  $$ T = \left(\matrix {z^{j_1} + ua    &  uc
 \cr ud  &  z^{j_2} + ub  \cr }\right)$$
where a, b, c, and d are holomorphic functions in
$U \cap V.$
We will change  coordinates to obtain an upper triangular
 transition matrix 
$$ \left(\matrix {z^{j_1}   & uc  \cr 0  &  z^{j_2}   \cr }\right),$$
which is equivalent to an extension
$$0 \rightarrow {\cal O}^!(-j_1) \rightarrow E \rightarrow {\cal O}^!(-j_2) \rightarrow 0.$$
Our required change of coordinates 
will be 
$$\left(\matrix {1  & 0  \cr \eta  &  1  \cr }\right)
\left(\matrix {z^{j_1} + ua  & uc  \cr ud  &  z^{j_2} + ub  \cr }\right)
\left(\matrix {1  & 0  \cr \xi  &  1  \cr }\right)$$
where $\xi$ is a holomorphic function on $U$
and $\eta$ is a holomorphic function on $V$
whose values 
 will be determined
in the following calculations.

After performing this multiplication, the entry $e(2,1)$ of
the resulting matrix   is
$$ e(2,1)= \eta\, (z^{j_1} + ua) +ud+ [\eta uc+(z^{j_2}+ ub
)]\,\xi.$$
We will
choose $\xi$ and $\eta $ to make  $e(2,1) = 0.$
We write the power series expansions for 
 $\xi$ and $\eta $ as 
$\xi = \sum_{i = 0} ^\infty \xi_i(z)\, u^i$
and $\eta = \sum_{i = 0} ^\infty \eta_i(z^{-1})\,(z^k u)^i,$
and plug into the expression for $e(2,1) .$ 
The term independent of $u$ in $e(2,1) $ is
$$\eta_0(z^{-1})z^{j_1} + \xi_0(z)z^{j_2}.$$
Since $j_2 - j_1 \le 0$ we may choose
$\eta(z^{-1}) = z^{j_2 - j_1}$ and  $\xi(z) =  1.$
After these choices $e(2,1)$ is now a multiple of $u.$
Suppose that the coefficients of $\eta$ and $\xi$
have been chosen up to power $u^{n-1}$ so that
$e(2,1)$ becomes a multiple of $u^n.$
Then  the coefficient of  $u^n$
in the expression for
$e(2,1)$ is
$$ \eta_n\,z^{j_1+kn} + \xi_n\,z^{j_2}+ \Phi$$
where $$\Phi =
\sum_{s+i = n}\eta_s\,a_i\,d_n\,z^{sk}+
\sum_{s+i+m = n}  \eta_s\,c_i\, \xi_m\,z^{sk}
 \sum_{m+i = n}\xi_m\,b_i.$$
We separate $\Phi$ into two parts
$$ \Phi = \Phi_{>j_2}   + \Phi_{\le j_2}$$
where $ \Phi_{>j_2}$ is the part of $\Phi $ containing
the powers $z^i $ for $i> j_2$ and
 $ \Phi_{\le j_2} $ is the part of $\Phi $ containing
powers  $z^i$ for $i \le j_2.$
We then choose the values of $\eta_n$ and $\xi_n$
as
$$\eta_n = z^{-j^1-nk} \Phi_{\le j_2} $$ and
$$\xi_n = z^{-j_2} \Phi_{>j_2} .$$
These choices cancel the coefficient of $u^n$ in
$e(2,1).$
Induction on $n$ gives $e(2,1)= 0$
And provides a  transition matrix of the form
 $$ T = \left(\matrix {z^{j_1} + ua    &  uc
 \cr 0 &  z^{j_2} + ub  \cr }\right).$$
Now do a similar trick using the  change of coordinates
 $$  \left(\matrix {\eta_1    &  0
                    \cr 0 &  \eta_2 \cr }\right)
     \left(\matrix {z^{j_1} + ua    &  uc
                    \cr 0 &  z^{j_2} + ub  \cr }\right)
      \left(\matrix {\xi_1    &  0
                     \cr 0 &  \xi_2 \cr }\right)$$
and choose $\xi_1,\, \xi_2,\, \eta_1$ and $\eta_2$ appropriately
to obtain a transition matrix of the form
 $$ T = \left(\matrix {z^{j_1}    &  uc
 \cr 0 &  z^{j_2}  \cr }\right).$$ \hfill\ja

\begin{theorem}:
Holomorphic bundles over $M_k, \,\, k>0$ are algebraic.
\end{theorem}

\noindent{\bf Proof}:
Let $E$ be a holomorphic bundle over
$M_k$ whose restriction to the exceptional 
divisor is 
$E_{\ell} \simeq  \oplus_{i = 1}^n {\cal O}(j_i), $
then  $E$ has a 
transition matrix
of the form
$$\left(\matrix {z^{j_1} & p_{12} & p_{13} & \cdots  
                \cr 0 &  z^{j_2} & p_{23} & p_{24} &  \cdots
                \cr \vdots &     & \vdots &     
                \cr 0  & \cdots & 0 & z^{j_{n-1}} &  p_{n-1,n}
                \cr 0 & \cdots & & 0 & z^{j_n} \cr}\right)$$
from $U$ to $V,$  where $p_{ij}$ are polinomials 
defined on $U \cap V.$

 Once again we  will give the detailed proof for the case 
$n = 2.$ The general proof is by induction on $n$ and is essentially 
the same as for $n = 2$ only notationally uglier.\hfill\ja

\vspace{5 mm}

For the case $n = 2$ we restate the theorem giving the specific form of the polinomial.

 \begin{theorem}:   Let $E$ be a
 holomorphic  rank two vector  bundle on
 $M_k$ whose restriction to the exceptional divisor is  
$E_{\ell} \simeq  {\cal O}(j_1)  \oplus {\cal O}(j_2),$ with $j_1 \ge j_2.$   
Then  $E$  has a transition matrix
of the form

$$\left(\matrix {z^{j_1} & p  
                \cr   0 & z^{j_2}  
                \cr 
}\right)$$
from $U$ to $V,$  where the polinomial $p$ is given by

$$p = \sum_{i = 1}^{ \left[(j_1 - j_2 -2)/k\right]}
 \sum_{l = ki+j_2+1}^{j_1-1}p_{il}z^lu^i$$ and 
$p = 0$ if $j_1< j_2 +2.$

\end{theorem}

\noindent {\bf Proof}:
Based on the proof of Theorem 2.1 we know that $E$
has a transition matrix of the form
$$  \left(\matrix {z^{j_1}    &  uc
                    \cr 0 &  z^{j_2}  \cr }\right).$$
We are left with obtaining the form of the polinomial $p,$
for which we perform another set of coordinate changes as follows.

 $$  \left(\matrix {1    &  \eta  
                    \cr 0 &  1  \cr }\right)
     \left(\matrix {z^{j_1}    &  uc  
                    \cr 0 &  z^{j_2}  \cr }\right)
     \left(\matrix {1    &  \xi  
                    \cr 0 &  1  \cr }\right),$$
where the coefficients of 
$\xi = \sum_{i = 0} ^\infty \xi_i(z)\, u^i$
and $\eta = \sum_{i = 0} ^\infty \eta_i(z^{-1})\,(z^k u)^i,$
will be choosen apropriately in  
the following steps.
After performing this multiplication, the entry
$e(1,2)$  of the resulting matrix is 
$$e(1,2) =  z^{j_1}\,\xi +uc + z^{j_2}\,\eta.$$
The term independent of $u$ in the expression
for  $e(1,2)$ is $z^{j_1}\xi_0(z) + z^{j_2+k} \eta_0(z^{-1}).$
However, we know from the expression
for our matrix $T$ (proof of lemma 2.1), that  $e(1,2)$
must be a multiple of $u;$ accordingly
we  choose 
$\xi_0(z) = \eta_0(z^{-1}) = 0.$
Placing this information into the above equation, we obtain
$$e(1,2) = \sum_{n = 1}^\infty (\xi_n(z)z^{j_1} +
c_n(z,z^{-1}) + \eta_n(z^{-1})z^ {j_2+kn})\,u^n.$$
Proceeding as we did in the proof of Lemma 
2.1, we choose values of $\xi_n$ and $\eta_n$
to cancel as many coefficients 
of $z$ and $z^{-1}$ as possible.
However, here $\xi_n$ appears multiplied by $z^{j_1}$
( and $\eta_n$ multiplied by $z^{j_2 + kn}$),
therefore the optimal choice of coeficients
 cancels only powers of $z^i$
with $i \ge j_1$ (resp. $z^i$ with
$i \le j_2 + kn$).
Consequently, $e(1,2)$ is left  only with
terms in $z^l$
for $j_2+nk < l< j_1 $, and we have the
expression
$$e(1,2) =  \sum_{i = 1}^\infty\sum_{l = nk+j_2+1}^{j_1-1}c_{il}z^lu^i.$$
But $i$ may only vary up to the point where
$nk+j_2+1 \leq j_1-1$
and the polynomial $p$ is given by
$$p = \sum_{i = 1}^{[(j_1-j_2-2)/k]}
 \sum_{l = ik+j_2+1}^{j_1-1}p_{il}z^lu^i.$$\hfill \ja

\section{Triviality outside the zero section}

From the previous section we know that bundles on $M_k$ are 
extensions of line bundles.  First we have the following lemma.

\begin{lemma}: Line bundles on $M_k$ are trivial outside the zero section.
\end{lemma}

\noindent{\bf Proof}:
A line bundle on $M_k$ can be given by a transition function $z^j$ for some 
integer $j.$  Then the function given by $z^{k-j}u$ on $U$ and 
$z^ku$ on $V$ is a global holomorphic section 
which trivializes the bundle
outside the zero section.
\hfill\ja

We now show that the extensions given in Section 2 are trivial outside 
the zero section.

\begin{theorem} Holomorphic vector bundles on ${\cal O}(-1)$ are trivial outside 
the zero section.
\end{theorem}

\noindent {\bf Proof}:Let $E$ be a holomorphic 
bundle on  ${\cal O}(-1).$
According to the previous section we know that 
$E$ is algebraic. Call $F$ the restriction of 
$E$ to the complement of the zero section, i.e. $F =  E|_{\ell^c}.$
Let $\pi : {\cal O}(-1) \rightarrow {\bf C}^2$ be the blow up
map. Then $\pi_*(F)$ is an algebraic bundle
over ${\bf C}^2 - {0}$ and therefore it extends 
to a coherent sheaf ${\cal F}$ over ${\bf C}^2.$
Then ${\cal F} ^{**}$ is a reflexive sheaf 
and as such has singularity set of codimension
3 or more, hence in this case ${\cal F} ^{**}$
is locally free. Moreover, as a bundle on 
${\bf C}^2$ it must be holomorphically  trivial.
 But ${\cal F} ^{**}$ restricts to
  $\pi_*(F)$ on ${\bf C}^2 - {0} ,$ hence
  $\pi_*(F)$ is trivial and so is $F.$ 
\hfill\ja

\begin{corollary} Holomorphic bundles on the blow up
of a surface are trivial on a  neighborhood 
of the exceptional divisor minus the 
exceptional divisor.
\end{corollary}

\noindent{\bf Proof}: Apply Theorem 3.2 to 
${\widetilde {\bf C}^2} = {\cal O}(-1)$.\hfill\ja

\vspace {5 mm}

\noindent Elizabeth Gasparim \\
International Centre for Theoretical Physics \\
P.O. Box 586 \\
34100 Trieste, ITALIA\\
gasparim@ictp.trieste.it\\
\end{document}